\documentclass[12pt,leqno]{article}
\usepackage[francais]{babel}
\usepackage[latin1]{inputenc}
\usepackage[T1]{fontenc}
\usepackage{amsmath,amsfonts,amssymb,amsthm,amscd}
\usepackage{bm}
\date{}

\newcommand\N{\mathbb{N}}

\newcommand{\newsmile}{\smallsmile_{\hspace{-3mm}\displaystyle\smallfrown}}

\def\adots{\mathinner{\mkern1mu\raise1pt\vbox{\kern7pt\hbox{.}}
\mkern2mu\raise4pt\hbox{.}
\mkern2mu\raise7pt\hbox{.}\mkern1mu}}

\begin{document}

\title{\'Etude du graphe divisoriel 4}

\author{Pierre Mazet et Eric Saias}


\maketitle

\section{Introduction}

Appelons  chaîne-permutation du graphe divisoriel de  $\N^*$ (ou plus simplement chaîne-permutation) toute application bijective $f: \N^* \rightarrow \N^* $ telle que $f(n)$ est un diviseur ou un multiple de $f(n+1)$ pour tout entier $n\geq1$. Il est facile de construire des chaînes-permutations. Voici les premières valeurs de l'une d'entre elles: 
$$1-2-(2\times3)-3-(3\times4)-4-(4\times5)-5-(5\times7)-7-(7\times8)-8\cdots$$
On a ici $\limsup_{n \rightarrow +\infty}{f(n)/n^2}=1/4$. Le résultat ci-dessous montre que la convergence vers l'infini de $f(n)$ peut être beaucoup plus lente.

\vspace{4mm}

\noindent \textbf{Théorème}.

\textsl{
Il existe une constante $c_1$ et une chaîne-permutation $f$ telles que pour tout $n\geq2$, on a 
$$f(n) \leq c_1 n (\log n)^2 .$$}

\vspace{4mm}

Intéressons-nous maintenant à des permutations de $\N^*$ qui ne sont pas nécessairement des chaînes-permutations. Notons $[a,b]$ le plus petit commun multiple des entiers $a$ et $b$. En 1983, Erdös, Freud et Hegyvari ont montré qu'il existe une constante $c_2$ et une permutation $f$ de $\N^*$ telles que pour tout $n\geq3$, on a 
$$[f(n),f(n+1)]\leq n \exp\{c_2 \sqrt{\log n }  \log \log n\}$$
(voir le théorème 3 de \cite{EFH}). En modifiant légèrement la permutation choisie par Erdös, Freud et Hegyvari, Chen et Ji \cite{CJ} ont montré en 2011 que l'on peut remplacer la fonction à l'intérieur de l'accolade par $(2\sqrt{2}+o(1))\sqrt{\log n \log \log n}$. Il découle immédiatement de notre théorème une amélioration de leur résultat.

\vspace{4mm}

\noindent \textbf{Corollaire}.

\textsl{
Il existe une constante $c_3$ et une permutation $f$ de $\N^*$ telles que pour tout $n\geq2$, on a 
$$[f(n),f(n+1)]\leq c_3 n  (\log n)^2.$$}

\vspace{4mm}

Signalons que Tenenbaum \cite{TEN2} conjecture que l'on peut remplacer l'expression à droite du signe $\leq$ dans cette dernière formule par 
$n (\log n)^{1+o(1)}.$ L'exposant $1$ de $\log n$ serait alors optimal car on sait que (théorème 3 de \cite{AEDD}) pour toute permutation $f$ de $\N^*$, on a 
$$\limsup_{n \rightarrow +\infty}\frac{[f(n),f(n+1)]}{n \log n} > 0.$$

\vspace{4mm}

Pour établir le théorème, on utilise une construction de chaîne finie due à Tenenbaum (voir le paragraphe 4 de \cite{TEN2}) et un résultat relatif aux entiers à diviseurs denses (théorème 1 de \cite{EDD1}) (pour la définition de ces termes, voir le paragraphe suivant). Même si cela ne nous est pas utile ici, signalons au passage que le récent travail de Weingartner \cite{WEI} fournit un équivalent asymptotique de la fonction de comptage des entiers à diviseurs $y$-denses, uniforme en $y$.

\section{Notations}

Appelons chaîne toute application injective $f:\N^*\rightarrow \N^*$ telle que $f(n)$ est un diviseur ou un multiple de $f(n+1)$ pour tout entier $n\geq1$. Appelons aussi chaîne finie de longueur $l$ tout $l$-uplet $\mathcal{C} = a_1,a_2,\cdots,a_l$ d'entiers positifs deux à deux distincts et tels que $a_i$ est un diviseur ou un multiple de $a_{i+1}$ pour tout entier $i$ de l'intervalle $[1,l-1]$. On notera \ $\mathrm{longueur\ }(\mathcal{C}):= l.$

\vspace{2mm}

On note pour tout entier $n\geq2, P^+(n)$ (respectivement $P^-(n)$) le plus grand (resp. petit) facteur premier de n. On pose de plus $P^+(1) = 1$. 

\vspace{2mm}

On note 

$$
F(n) := \left\{
\begin{array}{cl}
1 & (n=1) \\ 
\mathrm{max} \{ dP^-(d) : d \vert n , d > 1\} & (n\geq2).
\end{array} 
\right.
$$

On dit qu'un entier est à diviseurs $y$-denses si $F(n)\leq yn.$ Cette dénomination provient de l'identité

$$ \frac{F(n)}{n} = max_{1\leq i < \tau (n)} \frac{d_{i+1}(n)}{d_{i}(n)}$$

{\parindent=0pt où} $1=d_{1}(n)< d_{2}(n) < \cdots < d_{\tau (n)}(n)=n $ désigne la suite croissante des diviseurs de n (voir le lemme 2.2 de \cite{TEN1}).

\section{Preuve du théorème}

Soit $p$ un nombre premier.\par

Au paragraphe 4 de \cite{TEN2}, Tenenbaum construit une famille de chaînes finies $\Gamma(x,p)$ pour tout réel $x \geq 2$. On s'intéresse ici au cas particulier $x=2p^2$ pour lequel on a les quatre propriétés suivantes.  \par
{\parindent=0pt $\Gamma(8,2)=1-2.$}\par
{\parindent=0pt Pour tout $p \geq 3, \Gamma(2p^2,p)=1-p-\cdots-2.$}\par
{\parindent=0pt Pour tout $p, \Gamma(2p^2,p)$ est constitué d'entiers $m$ sans facteur carré et tels que $m \leq 2p^2, P^+(m) \leq p.$ }\par
{\parindent=0pt Pour tout $p, \Gamma(2p^2,p)$ contient tous les entiers $m$ sans facteur carré et tels que $F(m) \leq p^2.$ }\par
{\parindent=0pt Cette dernière }propriété combinée avec le théorème 1 de \cite{EDD1} entraîne l'existence d'une constante $c>0$ telle que pour tout $p$
\begin{equation}
\mathrm{longueur\ } (\Gamma(2p^2,p)) \geq c p^2/ \log p .
\end{equation}

Nous modifions légèrement $\Gamma(2p^2,p)$ en considérant la chaîne finie $\mathcal{D}(p)$ obtenue en déplaçant dans $\Gamma(2p^2,p)$ le $1$ initial pour le mettre en final, puis en multipliant le tout par le nombre premier $p^*$ suivant immédiatement $p$. On a donc pour $p\geq 3, \mathcal{D}(p)=p^*p-\cdots-2p^*-p^*.$

\vspace{2mm}
 
Soit $f_0$ la suite commençant à $1$, puis obtenue par concaténation des $\mathcal{D}(p)$.
$$f_0: 1-\mathcal{D}(2)-\mathcal{D}(3)-\mathcal{D}(5)-\mathcal{D}(7)-\cdots$$
On vérifie que $f_0$ est une chaîne formée d'entiers sans facteur carré. On pose $q_0=1$.

On construit alors par récurrence une suite croissante $(q_k)_{k\geq0}$ de nombres égaux à $1$ ou premiers et des chaînes $f_k$ construites à partir de $f_0$ en remplaçant, pour $p\leq q_k$, la chaîne finie $\mathcal{D}(p)$ par une chaîne finie $\mathcal{C}(p)$ à définir.

\vspace{2mm}

Supposons construits $q_0\leq q_1 \leq \cdots \leq q_{k-1}$, $\mathcal{C}(p)$ pour $p\leq q_{k-1}$ ainsi donc que $f_0,f_1,\cdots,f_{k-1}$.
 Si $k$ est dans l'image de $f_{k-1}$, on pose $q_{k}=q_{k-1}$, il n'y a pas de nouvelle chaîne finie $\mathcal{C}(p)$ à définir et donc $f_k=f_{k-1}$.

On suppose dorénavant que $k$ n'est pas dans l'image de $f_{k-1}$.

\vspace{2mm}

On choisit $q_k$ un nombre premier tel que  $q_k>q_{k-1}$ et 
\begin{equation}
q_k\geq k^2.
\end{equation}

On définit alors $\mathcal{C}(p)$ pour $q_{k-1} < p \leq q_k $ de la manière suivante.\par 
{\parindent=0pt Si $q_{k-1} < p < q_k$}, on pose $\mathcal{C}(p)=\mathcal{D}(p).$\par
{\parindent=0pt Si  $p= q_k$}, on choisit 
$$\mathcal{C}(p)=pk^2-k-p^*pk^2-\mathcal{D}(p).$$
Le fait que $\mathcal{D}(p)$ est formée d'entiers sans facteur carré assure que les éléments de $\mathcal{C}(p)$ sont deux à deux distincts et donc que $\mathcal{C}(p)$ est une chaîne finie.

\vspace{2mm}

On pose alors

$$f_k: 1-\mathcal{C}(2)-\mathcal{C}(3)-\cdots -\mathcal{C}(q_k)-\mathcal{D}(q_k^*)-\mathcal{D}(q_k^{**})-\cdots$$

{\parindent=0pt et on vérifie que c'est une chaîne. On définit enfin}

$$f=\lim_{k \rightarrow +\infty} f_k = 1-\mathcal{C}(2)-\mathcal{C}(3)-\mathcal{C}(5)-\cdots-\mathcal{C}(p)-\cdots$$

{\parindent=0pt et on vérifie que c'est une chaîne-permutation.}

\vspace{2mm}

Le début de cette chaîne-permutation est formée des trois entiers $1-6-3$, et est suivie de la chaîne finie $\mathcal{C}(3).$ Soit $n \geq 4;$ il existe donc un nombre premier $p$ tel que $f(n) \in \mathcal{C}(p^*)$. En utilisant notamment $(3.2)$, on a \par 
{\parindent=0pt $f(n) \leq 2p^{*2}p^{**}\newsmile p^3.$ }Par ailleurs notons $r$ un nombre premier générique. On a en utilisant (3.1) et le théorème des nombres premiers 
$$n\geq \sum_{r\leq p}\mathrm{longueur\ }(\mathcal{C}(r)) \gg \sum_{r\leq p} r^2/\log r\newsmile p^3/(\log p)^2.$$
On en déduit que $f(n) \ll n (\log n)^2$. Cela conclut la démonstration du théorème.

\vskip4mm

\begin{tabular}{ll}

Pierre Mazet & Eric Saias \\

 &Laboratoire de Probabilités, Statistique et Modélisation,\\

 &Sorbonne Université,\\

 &4, place Jussieu, 75252 Paris Cedex 05 (France)\\

\vspace{2mm}

\textsf{piermazet@laposte.net} &\textsf{eric.saias@upmc.fr}
\end{tabular}


\vspace*{2cm}



\begin{thebibliography}{6}

\bibitem{CJ} Y.-G.    \textsc{Chen} and C.-S. \textsc{Ji}. --- \textit{The permutation of integers with small least common multiple of two subsequent terms,} Acta Math. Hungarica 132 (2011), 307-309.

\bibitem{EFH} P.  \textsc{Erdös,} R. \textsc{Freud} and N.\textsc{Hegyvari}.--- \textit{Arithmetical properties of permutation of integers,} Acta Math. Hungarica 41 (1983), 169-176.

\bibitem{EDD1}   E. \textsc{Saias}.--- \textit{Entiers à diviseurs denses 1,} J. Number Theory 62 (1997), 163-191.

\bibitem{AEDD}   E. \textsc{Saias}.--- \textit{Applications des entiers à diviseurs denses,} Acta Arithmetica 83 (1998), 225-240.

\bibitem{TEN1} G.  \textsc{Tenenbaum}.--- \textit{Sur un problème de crible et ses applications,} Ann. Sci. Ecole Norm. Sup. (4) 19 (1986), 1-30.

\bibitem{TEN2} G.  \textsc{Tenenbaum}.--- \textit{Sur un problème de crible et ses applications, 2. Corrigendum et étude du graphe divisoriel,} Ann. Sci. Ecole Norm. Sup. (4) 28 (1995), 115-127.

\bibitem{WEI}   A. \textsc{Weingartner}.--- \textit{Practical numbers and the distribution of divisors,} Quart. J. Maths 66 (2015), 743-758.



\end{thebibliography}
\end{document}